\documentclass[dvips,twoside]{article}
\usepackage[latin1]{inputenc}
\usepackage[T1]{fontenc}
\usepackage{parskip}
\usepackage{amsmath,amsfonts,amssymb,amsxtra}
\usepackage{latexsym}
\usepackage{theorem}
\usepackage{graphicx,psfrag}
\usepackage{verbatim}

{\theoremstyle{plain}                       
\theorembodyfont{\itshape}
\newtheorem{Theorem}{Theorem}}

{\theoremstyle{plain}                       
\theorembodyfont{\itshape}
}

{\theoremstyle{plain}                       
\theorembodyfont{\itshape}
\newtheorem{Proposition}{Proposition}}

{\theoremstyle{plain}                       
\theorembodyfont{\rmfamily}
\newtheorem{Definition}{Definition}}

{\theoremstyle{plain}                       
\theorembodyfont{\rmfamily}
}

{\theoremstyle{plain}                       
\theorembodyfont{\rmfamily}
\newtheorem{Example}{Example}}

{\theoremstyle{plain}                       
\theorembodyfont{\rmfamily}
\newtheorem{Remark}{Remark}}

{\theoremstyle{plain}                       
\theorembodyfont{\itshape}
\newtheorem{Lemma}{Lemma}}

{\theoremstyle{plain}                       
\theorembodyfont{\itshape}
\newtheorem{Claim}{Claim}}

{\theoremstyle{plain}                       
\theorembodyfont{\rmfamily}
\newtheorem{Conjecture}{Conjecture}}

\begin{document}

\title{ Fundamental Markov systems}
\author{Ivan Werner\\ {\small Moscow, Russia}\\
   {\small Email: ivan\_werner@mail.ru}}
\maketitle

\begin{abstract}\noindent
 We continue development of the theory of Markov systems initiated
 in \cite{Wer1}. In this paper, we introduce fundamental Markov systems associated with random dynamical
 systems and show that the proof of the uniqueness and empiricalness
 of the stationary initial distribution of the random dynamical
 system reduces to that for the fundamental Markov system associated
 with it. The stability criteria for the latter are much
 clearer.


\end{abstract}

In \cite{Wer1}, the author initiated the study of a general concept
of a Markov system. This was motivated by a desire to have a
scientifically consistent unifying mathematical structure which
would cover {\it finite Markov chains} \cite{B}, $g$-measures
\cite{Le} and {\it iterated function systems with place-dependent
probabilities} \cite{BDEG}, \cite{El}.

The purpose of this note is to show that the structure of a Markov
system arises naturally (possibly unavoidably) in the study of {\it
random dynamical systems}.

\section{Random dynamical systems}

Let $(K,d)$ be a complete separable metric space and $E$ a countable
set. For each $e\in E$ let a  Borel measurable map
$w_e:K\longrightarrow K$ and a Borel measurable probability function
$p_e:K\longrightarrow [0,1]$ be given, i.e.
\[\sum_{e\in E}p_e(x)=1\;\;\;\mbox{ for all }x\in K.\] We call the
family $\mathcal{D}:=(K, w_e, p_e)_{e\in E}$ a {\it random dynamical system}. A
survey on random dynamical systems can be found e.g. in \cite{Kif}. If a reader
doesn't see, how the definition of random dynamical system in \cite{Kif} relates to
that in this paper, it is explained in \cite{KS}.

With the random dynamical system $\mathcal{D}$ is associated a {\it
Markov operator} $U$ acting on all bounded Borel measurable
functions $f$ by
\[Uf:=\sum\limits_{e\in E}p_ef\circ w_e\] and its {\it adjoint
operator} acting on the set of Borel probability measures $\nu$ by
\[U^*\nu(f):=\int U(f)\ d\nu.\]
A measure $\mu$ is called {\it invariant} with respect to the random
dynamical systems if and only if
\[U^*\mu=\mu.\]

Let $\Sigma^+:=\{(\sigma_1,\sigma_2,...):\ \sigma_i\in E,\
i\in\mathbb{N}\}$ endowed with the product topology of discreet
topologies and $S:\Sigma^+\longrightarrow\Sigma^+$ be the left shift
map. For $x\in K$, let $P_x$ be the Borel probability measure on
$\Sigma^+$ given by
\[P_x( _1[e_1,...,e_n]):=p_{e_1}(x)p_{e_2}\circ
w_{e_1}(x)...p_{e_n}\circ w_{e_{n-1}}\circ ...\circ w_{e_1}(x),\]
for every cylinder set $_1[e_1,...,e_n]:=\{\sigma\in\Sigma^+:\
\sigma_1=e_1,...,\sigma_n=e_n\}$, which is associated with the
Markov process generated by the random dynamical system with the
Dirac initial distribution $\delta_x$.

\begin{Remark}\label{er}
    Note that each map $w_e$  needs to be
    defined only on a subset of $K$ where its probability function $p_e$ is greater than zero.
    In this case, one obtains a random dynamical system on $K$ by extending the maps on the whole space
    arbitrarily.
\end{Remark}

\subsection{Historical context}

First, I will give some historical roots of the theory, and then, I will list some
works which in my view form the historical context of this paper.

It was not clear until the beginning of the twentieth century whether the
independence of random variables is a necessary condition for the {\it low of large
numbers} and the {\it central limit theorem} to hold. The breakthrough was a work by
A. A. Markov \cite{M} in which he has extended {\it the low of large numbers} to
dependent random variables (he also extended the central limit theorem to such
processes, see \cite{BLW} for a nice account on Markov's work and life). Markov
restricted himself to processes where each random variable depends only on the
previous one. Such a process, in case of a discreet state space, is generated by a
transition matrix (or a directed graph with probability weights) and an initial
distribution. These processes, now known as {\it Markov chains}, found many
applications.

From the work of Markov naturally arises the question whether the low of large
numbers and other limit theorems hold true for a more general class of dependent
processes. The truly next class of processes can be only those where the dependence
of a random variable on the past is not restricted to any number of previous
variables. The study of such processes was initiated (and motivated by applications)
by O. Onicescu and G. Mihoc \cite{OM}. Remarkably, they found a way of constructing
such processes without giving infinitely many rules of dependence of a random
variable on the values of all previous, which is of course very important for
applications. They called such processes {\it les cha\^{i}nes \`{a} liaisons
compl\`{e}tes}. Their work gave rise to the theory of {\it dependence with compete
connections}, in which many limit theorems have been proved \cite{IG}. What I call
in this paper {\it random dynamical system} can be seen as a special case of this
theory. However, it must be said that every stationary process with values in $E$
can already be generated by such a random dynamical system. Let's illustrate the
dependence in the notation of this paper. If maps $w_e$ are contractions on a
complete metric space, the past is coded by the maps to a point in topological space
$K$ and the probability of a value from $E$ of the next random variable is then
obtained as a function evaluated at that point,
\[P\left(X_{1} = e_1|X_0 = e_0, X_{-1} = e_{-1}, ...\right) =
p_{e_1}\left(\lim\limits_{n\to-\infty} w_{e_0}\circ
w_{e_{-1}}\circ...\circ w_{e_{-n}}(x_0)\right)\] for all
$(...,e_{-1},e_0,e_{1},...)\in E^\mathbb{Z}$ and $x_0\in K$.

One could argue that the development so far was guided mainly by the internal
mathematical logic. However, the main reason for the development of the mathematical
language still is the striving of Homo Sapiens for a better description of the world
outside (even if many modern mathematical craftsmen have no interest in science at
all). The random dynamical systems such as in this paper arose naturally also as
{\it learning models} \cite{Kar}, \cite{Br}. Consider an intelligent object, e.g. a
rat. Its state of "intelligence" $x$ is assumed to take some value from
$[0,1]\subset\mathbb{R}$. The object is "asked some questions" and its responses are
measured as '0' ("wrong") or '1' ("right"). If it gives an answer '0', its
"intelligence" is moved to $w_0(x) := 1/2x$. Otherwise, its "intelligence" is moved
to $w_1(x) := 1/2x+1/2$. It is natural to assume that a response '0' can happen with
some probability $p_0(x)$ (depending on the current state of "intelligence") and a
response '1' with probability $p_1(x) := 1 - p_0(x)$. One assumes that the function
$p_0$ does not change with the time (if it did, it's likely that it would have some
time average, which only would matter \cite{BEH}). Furthermore, probability
functions $p_0$ and $p_1$ easily can be obtained empirically if sufficiently many
measurements are made. Many seminal works on Markov processes generated by random
dynamical systems have been motivated by learning models.

The striving for a classification of stationary random processes brought attention
to such systems once again. After Ornstein \cite{Or} had shown that Kolmogorov-Sinai
entropy is a complete invariant on Bernoulli shifts and not complete on Kolmogorov
automorphisms \cite{OS}, arose an interest in measure preserving transformations
which are Kolmogorov and not Bernoulli. Of course, a natural candidate for that is a
strongly mixing process on a finite alphabet with infinite memory, which is known to
be generated by these kind of random dynamical systems. It was M. Keane \cite{Ke}
who drew attention of ergodic theory people to such systems. He considered a special
case where the maps of the system can be obtained as inverse branches of a
(expanding) homomorphism $S$ of a compact metric space $K$ (e.g.
$K:=E^{\mathbb{N}\cup{0}}$ and $S$ is shift map, then $w_e(...,y_{-1},y_0) :=
(...,y_{-1},y_{0},e)$ for all $(...,y_{-1},y_0)\in K$). In this case, the
probabilities for these maps can be given in terms of a function
$g:K\longrightarrow\ [0,1]$ ($p_e(...,y_{-1},y_0) := g(...,y_{-1},y_{0},e)$ for all
$(...,y_{-1},y_0)\in K$) with the property that $\sum_{y\in S^{-1}(\{x\})}g(y)=1$
for all $x\in K$). Then the stationary process also is invariant with respect to
homomorphism $S$. He called the stationary process in this case a {\it $g$-measure}.
Much of the progress on the subject has been made from the study of this special
case. Of course, the one-sided symbolic space is more than just an example. However,
the restriction to this case clearly encourages an algebraic rather than a geometric
approach, which could be one of the reasons why the structure of {\it the
fundamental Markov system}, which is going to be introduced in this paper, can not
have been seen before.

A careful reader probably has noticed that in all of the examples so far the maps of
the system have been contractive. It was shown by J. Hutchinson \cite{Hut} that a
family of contractive maps on a compact metric space has a unique invariant subset.
These subsets are like footprints of this kind of "animals". They are often very
irregular on every scale. (This must be very fascinating already for itself because
there are even people who study these footprints without actually any interest in
any "animal".) Advances in computational technology allowed to produce beautiful
colorful pictures of such sets. Naturally, arose the idea of storing images with
such systems \cite{EY}. Moreover, contractive maps allow to code strings of symbols
such as '0' and '1' to a point in the state space. This can be useful for data
compression \cite{BDX}. This all attracted a new interest to the random dynamical
systems. The study of such systems for these purposes was continued mainly by M. F.
Barnsley {\it et al.} \cite{BDEG}, \cite{BDEGE}. They called them {\it iterated
function systems with place dependent probabilities}. Much of the literature on the
subject is now available also under this name.

Recently, such system have been suggested also as a natural model for a quantum
measurement process \cite{Sl}. Let $K$ be a state space of a quantum system and $E$
a set of possible outputs of a measurement apparatus. When the quantum system being
in a state $x\in K$ interacts with the measurement apparatus and we observe $e\in
E$, it is known that this usually results in the system moving to a new state $y\in
K$. Hence,  it is natural to assume that there is a map $w_e: K\rightarrow K$
associated with each output $e\in E$. Moreover, it is natural to assume that $e$
could be observed with some probability $p_e(x)$ in this experiment. This all also
can be formulated using the conventional language of Hilbert space and projection
operators \cite{Sl}. The reader probably has noticed the analogy with the setup of
the "learning model" above. It probably is an indication on the universality of the
approach. In fact, studying a system by asking it some "question" is a very natural
approach.

I hope I have convinced the reader that the structure which is being
studied in this paper is not one of the toys designed by some
mathematicians to keep them busy, but one which naturally has been
crystalizing in the modern science. It clearly needs to be properly
integrated into the body of mathematics which describes
deterministic, random and quantum paradigms. Now, I will list some
important contributions which were made for the purpose of
understanding such systems with respect to their stability and
ergodic properties which clearly form the historical context of this
paper.

W. Doeblin and R. Fortet \cite{DF} (1937) gave fairly weak condition
on strictly positive probability functions on a compact metric space
which insures that the system with contractive maps has a unique
(attractive) stationary initial distribution. In particular, this
condition is satisfied if the probabilities have a summable
variation (Dini-continuous).

L. Breiman \cite{Br} (1960) proved the strong low of large numbers
for Markov operators with the Feller property which posses a unique
stationary initial distribution on a compact Hausdorf space (see a
sharper result in \cite{Wer8}).

R. Isaac \cite{Is} (1962) introduced the average contractiveness
condition, which insures the uniqueness of the stationary state (he
proved it on a compact metric space with strictly positive
probability functions satisfying Lipschitz continuity).

F. Ledrappier \cite{Le} (1974) identified the $g$-measures as
projections of some equilibrium states defined by variational
principle with respect to the potential $\log g$ seen as a function
on $E^{\mathbb{Z}}$ (see \cite{Wer9} for the explanation in the
general case). Furthermore, he showed that the natural extension of
the $g$-measure is weakly Bernoulli if the $g$-function is strictly
positive and Dini-continuous. P. Walters \cite{Wal} (1975) extended
the result to such $g$ on a subshift of finite type (Markov system).

T. Kaijser  \cite{Kij} (1981) introduced  a local contractiveness on
average condition in the general setup of random systems with
complete connections, which also can imply the uniqueness of the
stationary state. He called his systems {\it weakly distance
diminishing random systems with complete connections}.

H. Berbee \cite{Ber} (1987) showed the uniqueness and the very weak Bernoulli
property of the $g$-measure for strictly positive $g$-functions on a full shift
satisfying a continuity condition which is weaker than the Dini-continuity. After
that, many other works have been devoted to the weakening of the algebraic
expression forcing a continuity of the $g$-function, but mostly only for the proof
of uniqueness of the stationary state (see  \"{O}. Stenflo \cite{Ste} (2003) and N.
Berger, Ch. Hoffman, V. Sidoravicius \cite{BHS} (2005) and the references there).
(It must be pointed out that the continuity of the probability functions is not
fundamental for the stability of such systems (e.g. Example \ref{E2})).

J. H. Elton \cite{El} (1987) recognized the importance of the
relation for $x,y\in K$ given by the equivalence of measures $P_x$
and $P_y$ for the proof of the ergodic theorem for such systems
(though, he still assumed the uniqueness of the stationary state,
which was shown later \cite{Wer7} to be not necessary). He proved
that this equivalence relation holds true for all $x,y\in K$ if all
probability functions are bounded away form zero, Dini-continuous
and the system satisfies a contractiveness on average condition (in
the language of this paper, the fundamental Markov system associated
with such a random dynamical system has a single vertex set).

M. F. Barnsley, S. G. Demko, J. H. Elton and J. S. Geronimo
\cite{BDEG} \cite{BDEGE}  (1989) made, in my view, two important
contributions to understanding the conditions for the stability of
such systems. They showed that the condition of the strict
positivity of the probability functions can be weakened (is not
fundamental). Secondly, they found implicitly a way of reduction of
the multiplicative average contractiveness condition to the additive
average contractiveness condition, though they did not accomplish it
completely (see \cite{Wer1} for details).

I. Werner \cite{Wer7} (2005) showed that the condition of equivalence of measures
$P_x$ and $P_y$ for all $x$ and $y$ in the same vertex set of a Markov system (see
next section), which is continuous, irreducible and contractive, is sufficient for
the uniqueness of the stationary initial distribution. For example (the example was
given by using some ideas of  A. Johansson and A. \"{O}berg \cite{JO}), this
condition is satisfied if the probability functions have a square summable variation
on each vertex set and are bounded away from zero, which is more general than
Elton's \cite{El} example.

This list is far from being complete. There are many other works,
which a reader can easily find under code names {\it $g$-measures,
iterated functions systems with place-dependent probabilities,
random systems with complete connections, random dynamical systems}
and {\it Markov systems}. A reader interested in the study of
general Markov operators is referred to \cite{Sz}.

\subsection{Markov systems}

Now, let us consider a special random dynamical system which we call
a {\it Markov system} \cite{Wer1}.

Let $K_1,K_2,...,K_N$ be a partition of a metric space $K$ into
non-empty Borel subsets  (we do not exclude the case $N=1$).
Furthermore, for each $i\in\{1,2,...,N\}$, let
\[w_{i1},w_{i2},...,w_{iL_i}:K_i\longrightarrow K\] be a family of Borel
measurable  maps such that for each $j\in\{1,2,...,L_i\}$ there
exists $n\in\{1,2,...,N\}$ such that $w_{ij}\left(K_i\right)\subset
K_n$ (Fig. 1). Finally, for each $i\in\{1,2,...,N\}$, let
\[p_{i1},p_{i2},...,p_{iL_i}:K_i\longrightarrow\mathbb{R}^+\]
be a family of positive Borel measurable probability functions
(associated with the maps), i.e. $p_{ij}>0$ for all $j$ and
$\sum_{j=1}^{L_i}p_{ij}(x)=1$ for all $x\in K_i$.

\begin{center}
\unitlength 1mm
\begin{picture}(70,70)\thicklines
\put(35,50){\circle{20}} \put(10,20){\framebox(15,15)}
\put(40,20){\line(2,3){10}} \put(40,20){\line(4,0){20}}
\put(50,35){\line(2,-3){10}} \put(5,15){$K_1$} \put(34,60){$K_2$}
\put(61,15){$K_3$} \put(31,50){\framebox(7.5,5)}
\put(33,45){\framebox(6.25,9.37)} \put(50,28){\circle{7.5}}
\put(45,21){\framebox(6,5)} \put(10,32.5){\line(6,1){15}}
\put(10,32.5){\line(3,-5){7.5}} \put(17.5,20){\line(1,2){7.5}}
\put(52,20){\line(2,3){4}} \put(13,44){$w_{e_1}$}
\put(35,38){$w_{e_2}$} \put(49,42){$w_{e_3}$}
\put(33,30.5){$w_{e_4}$} \put(30,15){$w_{e_5}$}
\put(65,37){$w_{e_6}$} \put(15,5){ Fig. 1. A Markov system.}
\put(0,60){$N=3$} \thinlines \linethickness{0.1mm}
\bezier{300}(17,37)(20,46)(32,52)
\bezier{50}(32,52)(30.5,51.7)(30,49.5)
\bezier{50}(32,52)(30,51)(28.7,51.7)
\bezier{300}(26,31)(35,36)(35,47)
\bezier{50}(35,47)(35,44.5)(33.5,44)
\bezier{50}(35,47)(35,44)(36,44) \bezier{300}(43,50)(49,42)(51,30)
\bezier{50}(51,30)(50.5,32)(49.2,32.6)
\bezier{50}(51,30)(50.6,32)(51.5,33.2)
\bezier{300}(39,20)(26,17)(18,25)
\bezier{50}(18,25)(19.5,24)(20,21.55)
\bezier{50}(18,25)(20,23.5)(22,24) \bezier{300}(26,26)(37,28)(47,24)
\bezier{50}(47,24)(45,25)(43,24) \bezier{50}(47,24)(45,25)(44,26.5)
\bezier{100}(54.5,31.9)(56,37.3)(61,36.9)
\bezier{100}(61,36.9)(64.5,36.5)(66,34)
\bezier{100}(66,34)(68,30.5)(64.9,26.8)
\bezier{100}(64.9,26.8)(61.6,23.3)(57,23)
\bezier{50}(57,23)(58.5,23.3)(60.1,22.7)
\bezier{50}(57,23)(58.8,23.3)(59.5,24.8)
\end{picture}
\end{center}

\begin{Definition}
 We call  $V:=\{1,...,N\}$ {\it the set of vertices} and the
 subsets $K_1,...,K_N$ are called {\it the vertex sets}. Further, we
 call
  \[E:=\left\{(i,n_i):i\in\{1,...,N\},n_i\in\{1,
...,L_i\}\right\}\]  {\it the set of edges} and we use the following
notations:
\[p_e:=p_{in}\mbox{ and }w_e:=w_{in}\mbox{ for }e:=(i,n)\in E.\]
Each edge is provided with a direction (an arrow) by marking
 {\it an initial vertex} through the map
 \begin{eqnarray*} i:E&\longrightarrow& V\\
                         (j,n)&\longmapsto& j.
 \end{eqnarray*}
  {\it The terminal vertex} $t(j,n)\in V$ of an edge $(j,n)\in E$ is
  determined by the corresponding map through
\begin{eqnarray*}
 t((j,n))&:=&k\ :\iff\ w_{jn}\left(K_j\right)\subset K_k.
\end{eqnarray*}
We call the quadruple $G:=(V,E,i,t)$  {\it a directed (multi)graph}
or {\it digraph}. A sequence (finite or infinite)
$(...,e_{-1},e_0,e_1,...)$ of edges which corresponds to a walk
along the arrows of the digraph (i.e. $t(e_k)=i(e_{k+1})$) is called
{\it a path}.
\end{Definition}
\begin{Definition}
 We call the family $\mathcal{M}:=\left(K_{i(e)}, w_e, p_e\right)_{e\in E}$ a (finite) {\it
  Markov system}.
 \end{Definition}
The definition can be easily generalized to the infinite case.

The Markov system defines a random dynamical system on $K$ by
extending the probability functions $p_e|_{K_{i(e)}}$ on the whole
space by zero and the maps arbitrarily, as in Remark \ref{er}.
 \begin{Definition} We call a Markov system  {\it irreducible} or
 {\it aperiodic} if and only if its directed graph is irreducible or
 aperiodic respectively.
\end{Definition}
\begin{Definition}[CMS]
We call Markov system $\mathcal{M}$ {\it contractive} with an {\it
average contracting rate} $0<a<1$ if and only if it satisfies the
following {\it condition of contractiveness on average}:
\begin{equation}\label{cc}
    \sum_{e\in E}p_{e}(x) d(w_{e}(x),w_{e}(y))\leq ad(x,y)\mbox{
for all } x,y\in K_{i}\mbox{ and }i\in\{1,...,N\}
\end{equation}
(it is understood here that $p_e$'s are extended on the whole space
by zero and $w_e$'s arbitrarily). This condition was discovered by
Richard Isaac in 1962 for the case $N=1$ \cite{Is}.
\end{Definition}

It was shown in \cite{Wer7} that an irreducible contractive Markov
system $\mathcal{M}$ with uniformly continuous probabilities
$p_e|_{K_{i(e)}}>0$ has a unique invariant Borel probability measure
if $P_x\ll P_y$ for all $x, y\in K_{i(e)}$, $e\in E$, and the
subsets $K_i$ form an open partition of $K$ (this was shown in
\cite{Wer7} for some locally compact spaces, but it holds also on
complete separable spaces, as contractive $\mathcal{M}$ also posses
invariant measures on such spaces \cite{HS}).

\subsection{Fundamental Markov systems}

Now, we intend to show that with every random dynamical system
$\mathcal{D}$ is associated an equivalent Markov system
$\mathcal{M}':=(K'_{i(e)}, w'_e, p'_e)_{e\in E'}$ (not necessarily
finite) such that  $P'_x\ll P'_y$ for all $x, y\in K'_{i(e)}$, $e\in
E'$, and each $K'_{i(e)}$ is the largest with such property, where
$P'_x$ are the probability measures on the code space of
$\mathcal{M}'$.

The construction of $\mathcal{M}'$ goes as follows. Define an
equivalence relation between $x,y\in K$ by
\[x\sim y\;\;\; :\Leftrightarrow\;\;\; P_x\ll\gg P_y,\]
where $P_x\ll\gg P_y$ means $P_x$ is absolutely continuous with
respect to $P_y$ and $P_y$ is absolutely continuous with respect to
$P_x$.  Let
\[\biguplus\limits_{i\in V'}K'_i=K\]
be the partition of $K$ into the equivalence classes.  Then, for
every $e\in E$ and $x,y\in K'_i$, $i\in V'$,
\[p_e(x)=0\;\;\;\Leftrightarrow\;\;\;p_e(y)=0.\]
Hence, for every $e\in E$ and $i\in V'$,
\begin{equation}\label{fp}
    \mbox{either }p_e|_{K'_i}=0\mbox{ or }p_e|_{K'_i}>0.
\end{equation}
 Furthermore, holds the following.
\begin{Proposition}\label{fms}
    For every $e\in E$ and $i\in V'$ with $p_e|_{K'_i}>0$, there
    exists $j\in V'$ such that $w_e(K'_i)\subset K'_j$.
\end{Proposition}
{Proof.} Let $x,y\in K'_i$. Observe that
\[P_x( _1[e,\sigma_1,...,\sigma_n])=p_e(x)P_{w_e(x)}(
_1[\sigma_1,...,\sigma_n])\] for every cylinder set
$_1[\sigma_1,...,\sigma_n]$. Hence,
\[P_{w_e(x)}(B)=\frac{P_x(S^{-1}(B)\cap _1[e])}{p_e(x)}\]
for every Borel $B\subset\Sigma^+$. Since the analogous formula holds true also for
$P_{w_e(y)}$, we conclude that $w_e(x)\sim w_e(y)$. Thus, there exists $j\in V'$
such that $w_e(K'_i)\subset K'_j$.\hfill$\Box$

By (\ref{fp}) and Proposition \ref{fms}, we can define a  Markov
systems associated with $\mathcal{D}$.
\begin{Definition}
 Let \[E'_i:=\{(i,e):\ p_e|_{K'_i}>0,\ e\in E\}\;\;\;\mbox{ for all }i\in V'\] and
 \[E':=\bigcup_{i\in V'}E'_i.\]
 For every $(i,e)\in E'$ set $p'_{(i,e)}:=p_e1_{K'_i}$, $w'_{(i,e)}:=w_e|_{K'_i}$, $i'((i,e))=i$ and $t'((i,e))=j$ where
 $w_e(K'_i)\subset K'_j$. Then
 $G':=(V',E',i',t')$ is a  directed graph and we call
 $\mathcal{M}':=(K'_{i(e)}, w'_e, p'_e)_{e\in E'}$  the {\it
 fundamental Markov systems} associated with the  random dynamical
 system $\mathcal{D}$.
\end{Definition}

Now, we need to show that the vertex sets of the fundamental Markov
system $\mathcal{M}'$ associated with $\mathcal{D}'$ are measurable.
Otherwise, possible Banach-Tarski effects might make our
construction scientifically irrelevant. For that, we need to make
clear the constructive nature of the equivalence relation which
defines the vertex sets.

For $x, y\in K$, let
\begin{equation*}
    X_n(\sigma):=\left\{\begin{array}{cc}\frac{P_x( _1[\sigma_1,...\sigma_n])}{P_y( _1[\sigma_1,...\sigma_n])}
    , & P_y( _1[\sigma_1,...\sigma_n])>0 \\
    0,& P_x( _1[\sigma_1,...\sigma_n])=0\\
    \infty,& P_x( _1[\sigma_1,...\sigma_n])>0\mbox{ and }P_y( _1[\sigma_1,...\sigma_n])=0
\end{array}\right.
\end{equation*}
and
\begin{equation*}
    Y_n(\sigma):=\left\{\begin{array}{cc}\frac{P_y( _1[\sigma_1,...\sigma_n])}{P_x( _1[\sigma_1,...\sigma_n])}
    , & P_x( _1[\sigma_1,...\sigma_n])>0 \\
    0,& P_y( _1[\sigma_1,...\sigma_n])=0\\
    \infty,& P_y( _1[\sigma_1,...\sigma_n])>0\mbox{ and }P_x( _1[\sigma_1,...\sigma_n])=0
\end{array}\right.
\end{equation*}
for all $\sigma\in\Sigma^+$. Define
\[\xi(x,y):=\limsup\limits_{M\to\infty}\sup\limits_{n\in\mathbb{N}}P_x(X_n>M)+
   \limsup\limits_{M\to\infty}\sup\limits_{n\in\mathbb{N}}P_y(Y_n>M).\]
Observe that each $x\longmapsto P_x( _1[\sigma_1,...\sigma_n])$ is a
Borel measurable function. Therefore, each $x\longmapsto P_x(X_n>M)$
is a Borel measurable function. Hence, $x\longmapsto \xi(x,y)$ is a
Borel measurable function for all $y\in K$. By the symmetry, also
$y\longmapsto \xi(x,y)$ is a Borel measurable function for all $x\in
K$.

\begin{Lemma}\label{el}
For all $x,y\in K$,
\[\xi(x,y)=0\;\;\;\mbox{ if and only if }\;\;\;P_x\ll\gg P_y.\]
\end{Lemma}
{\it Proof.} Let $x,y\in K$. Let $\mathcal{A}_n$ be the finite
$\sigma$-algebra on $\Sigma^+$ generated by the cylinders
$_1[\sigma_1,...\sigma_n]$. Now, observe that, for all $m\leq n$ and
$C_m\in\mathcal{A}_m$,
\begin{equation}\label{me}
\int\limits_{C_m} X_n\ dP_y=\sum\limits_{C_n\subset
C_m}P_x(C_n)=P_x(C_m)=\int\limits_{C_m}X_m\ dP_y.
\end{equation}
Hence, $(X_n,\mathcal{A}_n)_{n\in\mathbb{N}}$ is a $P_y$-martingale.
Analogously, $(Y_n,\mathcal{A}_n)_{n\in\mathbb{N}}$ is a
$P_x$-martingale. Moreover, by (\ref{me}),
\[P_x(X_n>M)=\int\limits_{\{X_n>M\}}X_n\ dP_y\]
and analogously
\[P_y(Y_n>M)=\int\limits_{\{Y_n>M\}}Y_n\ dP_x.\]
Hence,
\[\xi(x,y)=\limsup\limits_{M\to\infty}\sup\limits_{n\in\mathbb{N}}\int\limits_{\{X_n>M\}}X_n\ dP_y+
\limsup\limits_{M\to\infty}\sup\limits_{n\in\mathbb{N}}\int\limits_{\{Y_n>M\}}Y_n\
dP_x.\] Therefore,  $\xi(x,y)=0$ if and only if $X_n$ and $Y_n$ are
uniformly integrable martingales. Hence, the condition $\xi(x,y)=0$
implies that there exists $X\in\mathcal{L}^1(P_y)$ and
$Y\in\mathcal{L}^1(P_x)$ such that $X_n\to X$ and $Y_n\to Y$ both in
$\mathcal{L}^1$ sense, and $E_{P_y}(X|\mathcal{A}_m)=X_m$ $P_y$-a.e.
and $E_{P_x}(Y|\mathcal{A}_m)=Y_m$ $P_x$-a.e. for all $m$. Then, by
(\ref{me}),
\[\int\limits_{C_m} X\ dP_y=\int\limits_{C_m}X_m\ dP_y=P_x(C_m)\mbox{ for all }C_m\in\mathcal{A}_m.\]
Hence, the Borel probability measures $XP_y$ and $P_x$ agree on all
cylinder subsets of $\Sigma^+$, and therefore, are equal.
Analogously, $YP_x=P_y$. Thus, $P_x\ll\gg P_y$.

Conversely, $P_x\ll\gg P_y$ implies that $X_n$ and $Y_n$ are
uniformly integrable \cite{Bil}, i.e. $\xi(x,y)=0$.\hfill$\Box$

\begin{Remark}
 Note that it is not obvious from the definition of $\xi$ that the
 relation $\xi(x,y)=0$ is transitive.
\end{Remark}

\begin{Proposition}
 (i) The vertex sets $K'_i$, $i\in V'$, are Borel measurable.\\
 (ii) Consider all probability
 functions $p'_e|_{K'_{i(e)}}$, $e\in E'$, to be extended on $K$ by
 zero and all maps $w'_e|_{K'_{i(e)}}$, $e\in E'$, to be extended on
 $K$ arbitrarily.  Let $U'$ be the Markov operator associated with the Markov
system $\mathcal{M}'$. Then $U'=U$, i.e. $\mathcal{M}'$ is an
equivalent random dynamical system to $\mathcal{D}$.
\end{Proposition}
{\it Proof.} (i) Let $i\in V'$. Fix $y\in K'_i$ and set
$f(x):=\xi(x,y)$ for all $x\in K$. Then, by Lemma \ref{el},
$K'_i=f^{-1}(\{0\})$. Hence, as $f$ is Borel measurable, $K'_i$ is
Borel measurable.

(ii) Let $g$ be a bounded Borel measurable function on $K$ and $x\in
K$. Then there exists a unique $i\in V'$ such that $x\in K'_i$.
Hence, by the definition of $\mathcal{M}$,
\begin{eqnarray*}
U'g(x)&=&\sum\limits_{e\in E'}p'_e(x)g\circ
w'_e(x)=\sum\limits_{e\in E'_i}p'_e(x)g\circ
w'_e(x)=\sum\limits_{e\in E}p_e(x)g\circ w_e(x)\\
&=& Ug(x).
\end{eqnarray*}\hfill$\Box$

\begin{Example}
    Suppose the random dynamical system $\mathcal{D}$ is given by
    the contractive Markov system $\mathcal{M}$ such that the vertex sets $K_1,...,K_N$ form
    an open partition of $K$ and the probability functions $p_e|_{K_{i(e)}}$
    are bounded away from zero and have a square summable variation,
    i.e. $\sum_{n\in\mathbb{N}}\phi^2(a^n)<\infty$, where $\phi$ is
    the maximum of modules of uniform continuity of functions
    $p_e|_{K_{i(e)}}$, $e\in E$. Then, by Lemma 2 in \cite{Wer7},
    $P_x\ll P_y$ for all $x, y\in K_i$, $i=1,...,N$ (note that the openness of the partition was required in \cite{Wer7} only
    to insure that $\mathcal{M}$ has an invariant measure (Feller
    property)). Therefore, the fundamental Markov system associated
    with $\mathcal{M}$ is  $\mathcal{M}$ itself.
\end{Example}

\begin{Example}\label{E2}
    Let $\mathcal{D}_2 := ([0,1], w_e, p_e)_{e = 0,1}$ be the random dynamical system
    where $w_0(x) = x/3$, $w_1(x) = x/3+1/3$ for all $x\in [0,1]$,
    \begin{equation*}
        p_0(x) = \left\{\begin{array}{cc}
          0, & 0\leq x\leq\frac{1}{9} \\
          b, & \frac{1}{9}< x\leq 1
        \end{array}\right.,
    \end{equation*}
    with $0 < b < 1$, and $p_1 = 1 - p_0$.
    \begin{Claim} \label{Cl}
       The vertex sets of the fundamental Markov system associated
       with $\mathcal{D}_2$ are $K_0 := [0,1/9]$, $K_1 := (1/9, 1/3]$ and $K_2 := (1/3, 1]$.
    \end{Claim}
    {\it Proof.}
       First, observe that partition $K_0\cup K_1\cup K_2$ makes a Markov system
       with constant probabilities from $\mathcal{D}_2$, where the transition matrix associated
       with it is
       \begin{equation*}
        A_2 := \left(\begin{array}{ccc}
          0  &0 & 1 \\
          b  &0 &1-b\\
          0  &b &1-b
        \end{array}\right).
    \end{equation*}
    Therefore,
       \[ p_{e_n}(w_{e_{n-1}}\circ...\circ w_{e_1}x) = p_{e_n}(w_{e_{n-1}}\circ...\circ
       w_{e_1}y)\] for all $x,y\in K_i$, $i = 0, 1,
       2$ and $e_1,...,e_n\in \{0,1\}$.
    Thus
       \[ P_x = P_y\mbox{ for all }x,y\in K_0, K_1, K_2.\]
    It remains to show that $K_0$, $K_1$ and $K_2$ are the largest with the
    property that $P_x<<>>P_y$ for all $x,y\in K_i$, $i = 0,1,2$.
    Let $x\in K_0$ and $y\not\in K_0$. Then $P_x( _1[0]) = 0$, but $P_y( _1[0]) =
    b>0$. Now, let $x\in K_1$ and $y\not\in K_1$. If $y\in K_0$, then $P_y( _1[0]) = 0$, but $P_x( _1[0]) =
    b$. Otherwise, if $y\in K_2$, $P_x( _1[00]) = 0$, but $P_y( _1[00])
    = b$. The claim follows.
    \hfill$\Box$

    Now, we can apply Theorem 2 in \cite{Wer2} to an equivalent fundamental Markov system on a
    disconnected set, the vertex sets of which are $\tilde K_0 := [0,1/9]$, $\tilde K_1 := [2/9, 1/3]$ and
    $\tilde K_2 := [2/3, 1]$. (Note that there is a missprint in \cite{Wer2} on page 471. It should be $\tilde A:= D^{-1}A^tD$.)
    By Theorem 2 in \cite{Wer2}, $\mathcal{D}_2$ has a unique invariant Borel probability
    measure $\mu_2$ with $\mu_2(K_0) = b^2/(1+b+b^2)$, $\mu_2(K_1) = b/(1+b+b^2)$
    and $\mu_2(K_2) = 1/(1+b+b^2)$ and the Markov chain associated
    with $\mathcal{D}_2$ is geometrically ergodic with a relative rate
    of convergence in Monge-Kantorovich metric less or equal to $\max\{1/3,b\}^{1/2}$.

    If we replace $p_0$ with
       \begin{equation*}
        p_0(x) = \left\{\begin{array}{cc}
          0, & 0\leq x\leq\frac{1}{27} \\
          b, & \frac{1}{27}< x\leq 1
        \end{array}\right.,
      \end{equation*}
    the same argumentation as in the proof of Claim \ref{Cl} shows that
    the fundamental Markov system associated with modified $\mathcal{D}_2$ has four vertex sets
     $K_0 := [0,1/27]$, $K_1 := (1/27,1/9]$, $K_2 := (1/9, 1/3]$ and $K_3 := (1/3,
    1]$ with constant probabilities on them given by the transition Matrix
    \begin{equation*}
    A_3 := \left(
      \begin{array}{cccc}
        0 & 0 & 0 & 1 \\
        b & 0 & 0 & 1-b \\
        0 & b & 0 & 1-b \\
        0 & 0 & b & 1-b \\
      \end{array}
    \right).
    \end{equation*}
    Analogously, applying Theorem 2 in \cite{Wer2} gives that the modified
    $\mathcal{D}_2$
    has an attractive invariant probability measure $\mu_3$ with the weights on vertex
    sets
    $\mu_3(K_0) = b^3/(1+b+b^2+b^3)$,
    $\mu_3(K_1) = b^2/(1+b+b^2+b^3)$, $\mu_3(K_2) = b/(1+b+b^2+b^3)$
    and $\mu_3(K_4) = 1/(1+b+b^2+b^3)$. (Somewhat surprisingly, the estimation of the relative rate of
    convergence by Theorem 2 in \cite{Wer2} remains the same (the absolute value of eigenvalues of such transition matrices
    remains the same and equals $b$). Note that in case of $p_0 = b$ the rate of convergence to the stationary
    state in the Monge-Kantorovich metric is not greater than $1/3$.)
\end{Example}

\begin{Example}\label{E3}
    Let $\mathcal{D}_3 := ([0,1], w_e, p_e)_{e = 0,1}$ be the random dynamical system
    where $w_0(x) = x/3$, $w_1(x) = x/3+1/3$ for all $x\in [0,1]$,
    \begin{equation*}
        p_0(x) = \left\{\begin{array}{cc}
          \frac{1}{4}, & 0\leq x\leq\frac{1}{2} \\
          \frac{1}{3}, & \frac{1}{2}< x\leq 1
        \end{array}\right.
    \end{equation*}
     and $p_1 = 1 - p_0$. As $p_{\sigma_2}(w_{\sigma_1}x) =
     p_{\sigma_2}(w_{\sigma_1}x)$ for all $x, y\in [0,1]$,
     $\sigma_1,\sigma_2\in\{0,1 \}$ and $p_0>0$, follows $P_x<<P_y$ for all $x,y\in
     [0,1]$. Thus the fundamental Markov system associated with $\mathcal{D}_3$ has
     the single vertex set.
\end{Example}

\begin{Example}\label{E4}
    Let $\mathcal{D}_4 := (\mathbb{R}, w_e, p_e)_{e = 0,1}$ be the random dynamical system
    where $w_0(x) = x/2$, $w_1(x) = x/2+1/2$ for all $x\in \mathbb{R}$,
    \begin{equation*}
        p_0(x) = \left\{\begin{array}{cc}
          \frac{1}{4}, & x\in\mathbb{Q} \\
          \frac{1}{3}, & x\in \mathbb{R}\setminus \mathbb{Q}
        \end{array}\right.
    \end{equation*}
     and $p_1 = 1 - p_0$. If $x$ is rational, then all images of $x$ under $w_0$ and $w_1$
     are rational. Therefore, $P_x$ is the Bernoulli measure generated with probabilities
     $\{1/4, 3/4\}$. Analogously, for irrational $y$, $P_y$ is the Bernoulli
     measure generated with probabilities $\{1/3, 2/3\}$. Since $P_x$ and $P_y$
     can not be absolutely continuous (absolutely continuous ergodic measures are equal), the fundamental Markov system associated with $\mathcal{D}_4$ has
     two vertex sets $K_0 := \mathbb{Q}$ and $K_1 := \mathbb{R}\setminus\mathbb{Q}$.
\end{Example}

Having obtained the well defined fundamental Markov system associated with
$\mathcal{D}$, we can define a Borel probability measure $P'_x$ on
$\Sigma'^+:=\{(\sigma_1,\sigma_2,...):\ \sigma_i\in E'\ \forall i\in\mathbb{N}\}$
(provided with the product topology of discreet topologies) by
\[P'_x( _1[e_1,...,e_n]):=p'_{e_1}(x)...p'_{e_n}(
w'_{e_{n-1}}\circ...\circ w'_{e_1}x)\] for all cylinder sets
$_1[e_1,...,e_n]\subset\Sigma'^+$, for every $x\in K$.

Now, we need to establish a relation between measures $P'_x$ and
$P_x$.

Since for each $e'\in E'$,  $w'_{e'}|_{K'_{i(e')}}$ is a restriction
of a unique $w_e$, we can define a map $\psi:\ E'\longrightarrow E$
by $\psi(e')$=e. This gives a Borel-Borel-measurable map $\Psi:\
\Sigma'^+\longrightarrow \Sigma^+$ by
$(\Psi(\sigma))_i:=\psi(\sigma_i)$ for all $i\in\mathbb{N}$ and for
all $\sigma\in\Sigma'^+$.

\begin{Lemma}\label{pl}
    For every $x\in K$, $P_x=\Psi(P'_x)$.
\end{Lemma}
{\it Proof.} Let $ _1[e_1,...,e_n]\subset\Sigma^+$ be a cylinder
set. Then, by the definition of $\Psi$, we can write
$\Psi^{-1}(_1[e_1,...,e_n])$ as a disjoint union of some cylinder
sets $ _1[e'_1,...,e'_n]\subset\Sigma'^+$, i.e.
\[\Psi^{-1}(_1[e_1,...,e_n])=\biguplus\limits_{(e'_1,...,e'_n),\psi(e'_i)=e_i}\ _1[e'_1,...,e'_n].\]
Therefore,
\[\Psi(P'_x)(_1[e_1,...,e_n])=\sum\limits_{(e'_1,...,e'_n),
\psi(e'_i)=e_i} p'_{e'_1}(x)...p'_{e'_n}\circ
w'_{e'_{n-1}}\circ...\circ w'_{e'_1}(x).\] Now, observe that, by the
definition of $p'_e$,
\[\sum\limits_{e',\psi(e')=e} p'_{e'}=p_{e}.\] This implies that
\[\Psi(P'_x)(_1[e_1,...,e_n])= p_{e_1}(x)...p_{e_n}\circ
w_{e_{n-1}}\circ...\circ w_{e_1}(x)=P_x(_1[e_1,...,e_n]).\] Thus,
the claim follows.\hfill$\Box$

\begin{Proposition}\label{acp}
  For all $x,y\in K$,
  \[P_x\ll P_y\;\;\;\Leftrightarrow\;\;\;P'_x\ll P'_y.\]
\end{Proposition}
{\it Proof.} $'\Leftarrow'$ Let $A\subset\Sigma^+$  Borel measurable
such that $P_y(A)=0$. Then, by Lemma \ref{pl},
$P'_y(\Psi^{-1}(A))=0$. Hence, by Lemma \ref{pl} and the hypothesis,
$P_x(A)=P'_x(\Psi^{-1}(A))=0$. Thus $P_x\ll P_y$.

$'\Rightarrow'$ First, observe that, by the construction of $M'$,
for every cylinder set $ _1[e'_1,...,e'_n]\subset\Sigma'^+$ with
$P'_y( _1[e'_1,...,e'_n])>0$, there exists a unique cylinder set $
_1[e_1,...,e_n]\subset\Sigma^+$ such that $P'_y(
_1[e'_1,...,e'_n])=P_y( _1[e_1,...,e_n])$ and $\Psi(
_1[e'_1,...,e'_n])= _1[e_1,...,e_n]$. Hence
\begin{equation}\label{pu}
P'_y( _1[e'_1,...,e'_n])=P_y( \Psi( _1[e'_1,...,e'_n])).
\end{equation}
Now, let $B\subset\Sigma'$ Borel measurable such that $P'_y(B)=0$.
Let $\epsilon>0$. By the hypothesis, there exists $\delta>0$ such
that
\begin{equation}\label{vu}
    P_y(C)<\delta\;\;\;\Rightarrow\;\;\; P_x(C)<\epsilon
\end{equation}
for all Borel measurable $C\subset\Sigma^+$. By the Borel regularity
of $P'_x$, there exists a countable family of cylinder sets
$C_k\subset\Sigma'^+$, $k\in\mathbb{N}$, such that
$B\subset\bigcup_{k} C_k$ and $P'_y(\bigcup_{k} C_k)<\delta$. Since
we can write every finite union $\bigcup_{k}^m C_k$ as a disjoint
union of cylinder sets $\bigcup_{k}^{n_m} \tilde C_k$,
\[\sum\limits_{k=1}^\infty P'_y(\tilde C)=
P'_y(\bigcup\limits_kC_k)<\delta.\] Hence, by (\ref{pu}) and the
Lemma \ref{pl},
\[P_y(\bigcup\limits_k\Psi(C_k))\leq \sum\limits_{k=1}^\infty P_y(\Psi(\tilde
C))=\sum\limits_{k=1}^\infty P'_y(\tilde C)<\delta.\] Therefore, by
(\ref{vu}) and Lemma \ref{pl}, \[P'_x(B)\leq
P'_x(\Psi^{-1}(\Psi(B)))=P_x(\Psi(B))\leq
P_x(\bigcup\limits_k\Psi(\tilde C_k))<\epsilon.\] Since $\epsilon$
was arbitrary, this completes the proof.\hfill$\Box$

\begin{Theorem} \label{mt}
    Suppose $\mathcal{D}$ is a random dynamical system with finitely
    many uniformly continuous probability functions $p_e$ and continuous maps $w_e$ on a complete separable
    metric space $K$. Suppose that the fundamental Markov system
    associated with $\mathcal{D}$ has finitely many vertices, is irreducible and contractive.
    Then\\
    (i) $\mathcal{D}$ has a unique invariant Borel probability measure
    $\mu$.\\
    (ii) For every $x\in K$,
   \[\frac{1}{n}\sum\limits_{k=o}^{n-1}f\circ{w_{\sigma_k}\circ...\circ w_{\sigma_1}(x)}\to
    \int f\ d\mu\mbox{ for
    $P_x$-a.e. }\sigma\in\Sigma^+\] for all bounded continuous functions $f$.
\end{Theorem}
{\it Proof.} Apply  Theorem 4 in \cite{Wer7} for the fundamental
Markov system associated with $\mathcal{D}$ with the following
justifications. In \cite{Wer7}, the vertex sets were required to
form an open partition of a state space in which sets of finite
diameter a relatively compact. This was to insure that the Markov
operator has the Feller property and an invariant Borel probability
measure. Here, the Feller property is already given by $\mathcal{D}$
and the existence of invariant measures for our fundamental Markov
system associated with $\mathcal{D}$ (on complete separable space)
was shown in \cite{HS}. Also, it was required  in \cite{Wer7} that
the probability functions $p_e|_{K_{i(e)}}$ shall be bounded away
from zero, but it was only required for the proof that $P_x\ll P_y$
for all $x$ and $y$ in the same vertex set (Lemma 2 in \cite{Wer7}).
The latter is given here by the construction of the fundamental
Markov system associated with $\mathcal{D}$ and Proposition
\ref{acp}.

After obtaining the result for the fundamental Markov system
associated with $\mathcal{D}$, deduce the result for $\mathcal{D}$
by Proposition \ref{acp}.\hfill$\Box$

\begin{Remark}
 Note that the fundamental Markov system $\mathcal{M}'$ associated with
 $\mathcal{D}$ is contractive if $\mathcal{D}$ is contractive.
\end{Remark}

\begin{Remark}
 Recall that there exist contractive random dynamical systems with strictly positive continuous probability
 functions which have more than one probability measure \cite{BK},\cite{BHS}. By Theorem \ref{mt} the
 fundamental Markov system associated with such a random dynamical system can not
 have a single vertex set.
\end{Remark}

\begin{Conjecture}\nonumber
I believe that fundamental Markov systems resolve the question of
the necessary and sufficient condition for the stability of such
random dynamical systems, which has been open already for more than
70 years, in the following way. {\it The random dynamical system has
a unique invariant Borel probability measure if and only if the
fundamental Markov system associated with it is recurrent (every
vertex of it is reached from any other by a finite path)}. Note that
a recurrent Markov system is necessarily countable (every vertex of
it can be coded by a finite path).
\end{Conjecture}

\end{document}